\documentclass[11pt,a4paper]{article}
\usepackage{a4wide}
\usepackage{amsmath,amssymb}
\usepackage[english]{babel}
\usepackage{graphicx}
\usepackage[numbers]{natbib}
\usepackage{theorem}
\usepackage{bm,setspace}
\usepackage{paralist,multirow}
\usepackage{ifpdf}
\usepackage{hyperref}
\hypersetup{
    bookmarksopen=false,
    bookmarksnumbered=true,
}
\ifpdf
  \hypersetup{colorlinks=true,linkcolor=black,urlcolor=black,citecolor=black,pdfpagemode=UseOutlines,plainpages=false,pdfpagelabels}
\else
  \hypersetup{colorlinks=false}
\fi

\theorembodyfont{\upshape}
\theoremheaderfont{\bfseries}

\newtheorem{theorem}{Theorem}[section]

\newtheorem{remark}[theorem]{Remark}
\newtheorem{proposition}[theorem]{Proposition}

\newenvironment{proof}{\par\textbf{Proof:}%\\
}{\hfill$\square$\par}
\numberwithin{equation}{section}

\newcommand{\LST}[1]{\widetilde{#1}}
\newcommand{\B}{\LST{B}}
\newcommand{\Bexh}{\LST{B}^\textrm{exh}}
\newcommand{\BP}{\LST{\textit{BP}}}
\newcommand{\Pexh}{P^\textrm{exh}}
\newcommand{\R}{\LST{R}}
\newcommand{\ee}{\textrm{e}}

\newcommand{\z}{\mathbf{z}}

\newcommand{\E}{\mathbb{E}}
\renewcommand{\P}{\mathbb{P}}

\newcommand{\sd}{\textrm{sd}}
\renewcommand{\O}{\mathcal{O}}

\newcommand{\ii}{i=1,\dots,N}

\newcommand{\Btilde}{\LST{B}^*}
\newcommand{\Ptilde}{\LST{P}^*}

\newcommand{\Rtilde}{\LST{R}^*}
\newcommand{\LB}{\widetilde{\textit{LB}}}
\newcommand{\LC}{\widetilde{\textit{LC}}}
\newcommand{\WC}{\widetilde{\textit{WC}}}
\renewcommand{\L}{\widetilde{L}}
\newcommand{\ub}[1]{\mathbf{B_{\bm{\mathit{#1}}}}}
\newcommand{\ubG}[1]{\mathbf{B^G_{\bm{\mathit{#1}}}}}
\newcommand{\ubexh}[1]{\mathbf{B'_{\bm{\mathit{#1}}}}}
\newcommand{\C}{\LST{C}}
\newcommand{\W}{\widetilde{W}}

\providecommand{\href}[2]{#2}

\title{Waiting times in queueing networks with a single shared server\footnote{The research was done in the framework of the BSIK/BRICKS project, the European Network of Excellence Euro-NF, and of the project ``Service Optimization and Quality'' (SeQual), funded by the Dutch agency SenterNovem. }}
\author{M.A.A. Boon\footnote{Eurandom and Department of Mathematics and Computer Science, Eindhoven University of Technology, P.O. Box 513, 5600MB Eindhoven, The Netherlands}\\\href{mailto:marko@win.tue.nl}{marko@win.tue.nl} \and R.D. van der Mei\ \footnote{Department of Mathematics, Section Stochastics, VU University, De Boelelaan 1081a, 1081HV Amsterdam, The Netherlands and Centre for Mathematics and Computer Science (CWI), 1098 SJ Amsterdam, The Netherlands}\\\href{mailto:mei@cwi.nl}{mei@cwi.nl} \and E.M.M. Winands \footnote{
University of Amsterdam, Korteweg-de Vries Institute for Mathematics, Science Park 904, 1098 XH  Amsterdam, The Netherlands}\\\href{mailto:e.m.m.winands@uva.nl}{e.m.m.winands@uva.nl}}

\date{November, 2012}

\begin{document}
\maketitle

\begin{abstract}
We study a queueing network with a single shared server that serves the queues in a cyclic order. External customers arrive at the queues according to independent Poisson processes. After completing service, a customer either leaves the system or is routed to another queue. This model is very generic and finds many applications in computer systems, communication networks, manufacturing systems, and robotics. Special cases of the introduced network include well-known polling models, tandem queues, systems with a waiting room, multi-stage models with parallel queues, and many others. A complicating factor of this model is that the internally rerouted customers do not arrive at the various queues according to a Poisson process, causing standard techniques to find waiting-time distributions to fail. In this paper we develop a new method to obtain exact expressions for the Laplace-Stieltjes transforms of the steady-state waiting-time distributions. This method can be applied to a wide variety of models which lacked an analysis of the waiting-time distribution until now.

\bigskip\noindent\textbf{Keywords:} queueing network,  waiting times, customer routing, shared server, polling

\bigskip\noindent\textbf{Mathematics Subject Classification:}  60K25, 90B22
\end{abstract}

\section{Introduction}\label{introduction}

In this paper we study a queueing network served by a single shared server that visits the queues in a cyclic order. Customers from the outside arrive at the queues according to independent Poisson processes, and the service time and switch-over time distributions are general. After receiving service at queue $i$, a customer is either routed to queue $j$ with probability $p_{i,j}$, or leaves the system with probability $p_{i,0}$. We consider systems with mixtures of gated and exhaustive service. This model can be seen as an extension of the standard polling model (in which customers always leave the system upon completion of their service) by customer routing. Yet another view is provided by the notion that the system is a Jackson network with a dedicated server for each queue with the additional complexity that only one server can be active in the network simultaneously. %\textit{The goal of the present paper is the derivation of the waiting-time distribution in a queueing network with a single shared server}.

The possibility of re-routing of customers further enhances the already-extensive modelling capabilities of polling models, which find applications in diverse areas such as computer systems, communication networks, logistics, flexible manufacturing systems, robotics systems, production systems and maintenance systems (see, for example, \cite{boonapplications2011,Grillo1,levy1,takagi3} for overviews). Applications of the introduced type of customer routing can be found in many of these areas. In this regard, we would like to mention a manufacturing system where products undergo service in a number of stages or in the context of rework \cite{Grasman1}, a Ferry based Wireless Local Area Network (FWLAN) in which nodes can communicate with each other or with the outer world via a message ferry \cite{kavitha}, a dynamic order picking system where the order picker drops off the picked items at the depot where sorting of the items is performed \cite{gongdekoster08}, and an internal mail delivery system where a clerk continuously makes rounds within the offices to pick up, sort and deliver mail \cite{sarkar}.

In the past many papers have been published on special cases of the current network. In some of these papers distributional results are derived as well; the techniques used do, however, not allow for extension to the general setting of the current paper.
Some special case configurations are standard polling systems \cite{takagi3}, tandem queues \cite{nair,taube}, multi-stage queueing models with parallel queues \cite{katayama}, feedback vacation queues \cite{boxmayechiali97, takine}, symmetric feedback polling systems \cite{takagifeedback,takine}, systems with a waiting room \cite{alineuts84,takacsfeedback77}, and many others. In conclusion, one can say that the present research can be seen as a unifying analysis of the waiting-time distribution for a wide variety of queueing models.

The main contribution of this paper is the derivation of waiting-time distributions in queueing networks with a single roving server via the development of a new method. For this model we derive the Laplace-Stieltjes transform of the waiting-time distribution of an arbitrary (internally rerouted, or external) customer. Due to this intrinsic complexity of the model, studies in the past were restricted to queue lengths and \textit{mean} delay figures (see \cite{boonvdmeiwinandsRovingPER2011,sarkar,sidi1,sidi2}). A complicating, yet interesting, factor is that the combined process of internal and external arrivals violates the classical assumption of Poisson (or even renewal) arrivals, implying that traditional methods are not applicable. The basic idea behind the new method is that we explicitly compute a priori all \emph{future} service requirements upon arrival of a newly arriving customer. In doing so the prerequisites of the distributional form of Little's Law are overcome.

An important feature of the newly developed technique is that it can be applied to a myriad of models which lacked an analysis of the waiting-time distribution until now. One could apply the framework (possibly after some minor modifications) to obtain distributional results in all of the aforementioned special cases of the studied system \cite{alineuts84,boxmayechiali97, katayama,nair,takacsfeedback77,takagifeedback,takagi3,takine,taube} but also, for example, in a closed network  \cite{altman2}, in an $M/G/1$ queue with permanent and transient customers \cite{boxmacohen91}, in a network with permanent and transient customers \cite{armonyyechiali99},
%in a queueing network with a single shared server and multiple customer classes at each queue \cite{hirayama},
%in an $M/G/1$ queueing model with gated random order of service \cite{resingrietman,rietmanresing},
or in a polling model with arrival rates that depend on the location of the server  \cite{boonsmartcustomers2010,smartcustomers}.
Although we study a continuous-time cyclic system with gated or exhaustive service in each queue, we may extend all results - without complicating the analysis - to discrete time, to periodic polling, to batch arrivals, or to systems with different branching-type service disciplines such as globally gated service.

The structure of the present paper is as follows. In Section \ref{modelsection}, we introduce the model and notation. Section \ref{waitingtimesectiongated} analyses the waiting-time distribution of an arbitrary customer for gated service. In Section \ref{exhaustiveservice} we study the system with mixtures of gated and exhaustive service. In the penultimate section, we present some examples which show the wide range of applicability of the studied model. The final section of this paper contains a brief discussion.

\section{Model description and notation}\label{modelsection}

We consider a queueing network consisting of $N\geq1$ infinite buffer queues $Q_1,\dots,Q_N$. External customers arrive at $Q_i$ according to a Poisson arrival process with rate $\lambda_i$, and have a generally distributed service requirement $B_i$ at $Q_i$, with mean value $b_i := \E[B_i]$. In general we denote the Laplace-Stieltjes Transform (LST) or Probability Generating Function (PGF) of a random variable $X$ with $\LST{X}(\cdot)$. The queues are served by a single server in cyclic order. Whenever the server switches from $Q_i$ to $Q_{i+1}$, a random switch-over time $R_i$ is incurred, with mean $r_i$. The cycle time $C_i$ is the time between successive moments when the server arrives at $Q_i$. The total switch-over time in a cycle is denoted by $R=\sum_{i=1}^N R_i$, and its first two moments are $r:=\E[R]$ and $r^{(2)}:=\E[R^2]$. Indices throughout the paper are modulo $N$, so $Q_{1-N}$ and $Q_{N+1}$ both refer to $Q_1$. All service times and switch-over times are mutually independent. This queueing network can be modelled as a \emph{polling system} with the specific feature that it allows for routing of the customers: upon completion of service at $Q_i$, a customer is either routed to $Q_j$ with probability $p_{i,j}$, or leaves the system with probability $p_{i,0}$. Note that $p_{i,0}$ should be greater than 0 for at least one queue, to make sure that customers can leave the system eventually. Moreover, note that $\sum_{j=0}^N~p_{i,j}=1$ for all $i$, and that the transition of a customer from $Q_i$ to $Q_j$ takes no time. Since we consider the gated and exhaustive service disciplines, the model under consideration has a branching structure, which is discussed in more detail by Foss \cite{fossBranching} in the context of queueing models, and by Resing \cite{resing93} more specifically in the context of polling systems.  The total arrival rate at $Q_i$ is denoted by $\gamma_i$, which is the unique solution of the following set of linear equations:
\begin{equation*}
\gamma_i = \lambda_i + \sum_{j=1}^N \gamma_j p_{j,i},\qquad\ii.
\end{equation*}
The offered load to $Q_i$ is $\rho_i:=\gamma_i b_i$ and the total load is $\rho:=\sum_{i=1}^N \rho_i$. We assume that the system is stable, which means that $\rho$ should be less than one (see \cite{sidi2}).
%xThe total service time $\Btotal_i$ of a customer is the total amount of service given during the presence
%of the customer in the network. Its first moment, denoted by
%$\btilde_i$, is uniquely determined by the following set of linear equations:
%\begin{eqnarray}
%\btilde_i &=& b_i + \sum_{j=1}^N \btilde_j p_{i,j},\qquad\qquad i=1,\dots,N.\label{btildei}
%%\btilde_i^{(2)} &=& b_i^{(2)} + 2b_i\sum_{j=1}^N \btilde_j p_{i,j}+ \sum_{j=1}^N \btilde_j^{(2)} p_{i,j}.
%\end{eqnarray}
%The LST of $\Btotal_i$ is not discussed in the present paper, but can be obtained by solving a similar set of equations.

\section{Gated service}\label{waitingtimesectiongated}

In the present section we study the waiting-time distribution of an arbitrary customer for a system in which each queue receives gated service, which means that only those customers present at the server's arrival at $Q_i$ will be served before the server switches to the next queue. We define the waiting time $W_i$ of an arbitrary customer in $Q_i$ as the time between his arrival at this queue and the moment at which his service starts. As far as waiting times are concerned, a customer that is routed to another queue, say $Q_j$, upon his service completion is regarded as a new customer with waiting time $W_j$. The waiting-time distribution is found by conditioning on the numbers of customers present in each queue at an arrival epoch. To this end, we study the joint queue-length distribution at several embedded epochs in Section~\ref{jointQLsubsection}. In Sections~\ref{cycleTimesubsection} and~\ref{waitingTimesubsection} we use these results to successively derive the cycle-time distribution and the waiting-time distributions of internally rerouted customers and external customers.

\subsection{The joint queue-length distributions}\label{jointQLsubsection}

Sidi et al. \cite{sidi2} derive the PGFs of the joint queue-length distributions in all $N$ queues at visit beginnings, visit completions, and at arbitrary points in time. In order to keep this manuscript self-contained, we briefly recapitulate their approach, as it forms the starting point of our novel method to find the waiting time LSTs. There is one important adaptation that we have to make, which will prove essential for finding waiting time LSTs. We consider not only the customers in all $N$ queues, but we distinguish between customers standing \emph{in front of} the gate and customers standing \emph{behind} the gate (meaning that they will be served in the next cycle). Hence, we introduce the $N+1$ dimensional vector $\z=(z_1, \dots, z_N, z_G)$. The element $z_i$, $i=1,\dots,N$, in this vector corresponds to customers in $Q_i$ standing in front of the gate. The element $z_G$ at position $N+1$ is only used during visit periods. During $V_j$, the visit period of $Q_j$, it corresponds to customers standing behind the gate in $Q_j$. This makes the analysis of systems with gated service slightly more involved than systems with exhaustive service (discussed in the next section). Before studying the joint queue-length distributions, we briefly introduce some convenient notation:
\begin{align*}
\Sigma(\z)&=\sum_{j=1}^N\lambda_j(1-z_j),\\
\Sigma_i(\z)&=\lambda_i(1-z_G)+\sum_{j\neq i}\lambda_j(1-z_j),\\
P_i(\z) &= p_{i,0}+p_{i,i}z_G+\sum_{j\neq i}p_{i,j}z_j.
\end{align*}

\paragraph{Visit beginnings and completions.}
A cycle consists of $N$ visit periods, $V_i$, each of which is followed by a switch-over time $R_i$, for $i=1,\dots,N$. A cycle $C_i$ starts with a visit to $Q_i$ and consists of the periods $V_i, R_i, V_{i+1},\dots,V_{i+N-1}$, $R_{i+N-1}$. Let $P$ denote any of these periods. We denote the joint queue length PGF at the \emph{beginning} of $P$ as $\LB^{(P)}(\z)$. The equivalent at the \emph{completion} of period $P$ is denoted by $\LC^{(P)}(\z)$. Since the gated service discipline is a so-called \emph{branching-type} service discipline (see \cite{resing93}), we can express each of these functions in terms of $\LB^{(V_i)}(\z)$, for any $i=1,\dots,N$. These relations, which are sometimes called \emph{laws of motion}, are given below.

\begin{subequations}
\begin{align}
\LC^{(V_i)}(\z) &= \LB^{(V_i)}\Big(z_1, \dots, z_{i-1},\B_i\big(\Sigma_i(\z)\big)P_i(\z), z_{i+1}, \dots, z_N, z_G\Big),\label{lawsofmotion1}\\
\LB^{(R_i)}(\z) &= \LC^{(V_i)}(z_1, \dots, z_N, z_i),\\
\LC^{(R_i)}(\z) &= \LB^{(R_i)}(\z)\R_i\Big(\Sigma(\z)\Big),\\
\LB^{(V_{i+1})}(\z) &= \LC^{(R_i)}(\z),\\
&\vdots\nonumber\\
\LB^{(V_{i+N})}(\z) &= \LC^{(R_{i+N-1})}(\z).\label{lawsofmotionN}
\end{align}
\end{subequations}
Note the subtle difference between $\LC^{(V_i)}(\z)$ and $\LB^{(R_i)}(\z)$, due to the fact that the gate in $Q_i$ is removed after the completion of $V_i$, causing type $G$ customers to become type $i$ customers. In steady-state we have that $\LB^{(V_{i+N})}(\z) = \LB^{(V_{i})}(\z)$, implying that we have obtained a recursive relation for $\LB^{(V_{i})}(\z)$. Resing \cite{resing93} shows how a clever definition of immigration and offspring generating functions can be used to find an explicit expression for $\LB^{(V_{i})}(\z)$. For reasons of compactness we refrain from doing so in the present paper. Instead we want to point out that the recursive relation obtained from \eqref{lawsofmotion1}-\eqref{lawsofmotionN} can be differentiated with respect to the variables $z_1, \dots, z_N, z_G$. The resulting set of equations, which are called the \emph{buffer occupancy equations} in the polling literature, can be used to compute the moments of the queue-length distributions at all visit beginnings and completions.

\paragraph{Service beginnings and completions.}
We denote the joint queue length PGF at \emph{service} beginnings and completions in $Q_j$ by respectively $\LB^{(B_j)}(\z)$ and $\LC^{(B_j)}(\z)$. Since a customer may be routed to another queue upon his service completion, we define $\LC^{(B_j)}(\z)$ as the PGF of the joint queue-length distribution right \emph{after} the tagged customer in $Q_j$ has received service (implying that he is no longer present in $Q_j$), but \emph{before} the moment that he may join another queue (even though these two epochs take place in a time span of length zero). Eisenberg \cite{eisenberg72} observed the following relation, albeit in a slightly different model:
\begin{equation}
\LB^{(V_i)}(\z) + \gamma_i\E[C]\LC^{(B_i)}(\z)P_i(\z) = \LC^{(V_i)}(\z) + \gamma_i\E[C]\LB^{(B_i)}(\z).
\label{eisenberg}
\end{equation}
Equation \eqref{eisenberg} is based on the observation that each visit beginning coincides with either a service beginning, or a visit completion (if no customer was present). Similarly, each visit completion coincides with either a visit beginning or a service completion. The long-run ratio between the number of visit beginnings/completions and service beginnings/completions in $Q_i$ is $\gamma_i\E[C]$, with $\E[C]=\E[C_i]=r/(1-\rho)$. The distribution of the cycle time is given in the next subsection.

Furthermore, Eisenberg observes the following simple relation between the joint queue-length distribution at service beginnings and completions:
\begin{equation}
\LC^{(B_i)}(\z) = \LB^{(B_i)}(\z)\B_i\big(\Sigma_i(\z)\big)/z_i.\label{servicebeginningscompletions}
\end{equation}
Substitution of \eqref{servicebeginningscompletions} in \eqref{eisenberg} gives an equation which can be solved to express $\LB^{(B_i)}(\z)$ in $\LB^{(V_i)}(\z)$ and $\LC^{(V_i)}(\z)$.

\paragraph{Arbitrary moments.}
The PGF of the joint queue-length distribution at arbitrary moments, denoted by $\L(\z)$, is found by conditioning on the period in the cycle during which the system is observed $(V_1, R_1, \dots, V_N, R_N)$.
\begin{equation}
\L(\z) = \frac{1}{\E[C]}\sum_{j=1}^N\left(\E[V_j]\L^{(V_j)}(\z)+r_j\L^{(R_j)}(\z)\right),\label{Lz}
\end{equation}
with $\E[V_j] = \rho_j\E[C]$. In \eqref{Lz} the functions $\L^{(V_j)}(\z)$ and $\L^{(R_j)}(\z)$ denote the PGFs of the joint queue-length distributions at an arbitrary moment during $V_j$ and $R_j$ respectively:
\begin{align}
\L^{(V_j)}(\z) &= \LB^{(B_j)}(\z)\frac{1-\B_j\big(\Sigma_j(\z)\big)}{b_j\Sigma_j(\z)},
\label{jointQLduringVj}\\
\L^{(R_j)}(\z) &= \LB^{(R_j)}(\z)\frac{1-\R_j\big(\Sigma(\z)\big)}{r_j\Sigma(\z)}.
\label{jointQLduringSj}
\end{align}
The interpretation of \eqref{jointQLduringVj} and \eqref{jointQLduringSj} is that the queue length vector at an arbitrary time point in $V_j$ or $R_j$ is the sum of those customers that were present at the beginning of that service/switch-over time, plus vector of the customers that have arrived during the elapsed part of the service/switch-over time. For more details about the joint queue length and workload distributions for general branching-type service disciplines (in the context of polling systems, but also applicable to our model) we refer to Boxma et al. \cite{boxmakellakosinski2011}.

\subsection{Cycle-time distributions}\label{cycleTimesubsection}

In the remainder of this paper we present new results for the model introduced in Section \ref{modelsection}. We start by analysing the distributions of the cycle times $C_i$, $i=1,\dots,N$. The idea behind the following analysis is to condition on the number of customers present in each queue at the beginning of $C_i$ (and, hence, of $V_i$). The cycle will consist of the service of all of these customers, plus all switch-over times $R_i, \dots, R_{i+N-1}$, plus the services of all customers that enter during these services and switch-over times \emph{and} will be served \emph{before} the next visit beginning to $Q_i$. The cycle time for polling systems without customer routing is discussed in Boxma et al. \cite{boxmafralixbruin08}. However, as it turns out, the analysis is severely complicated by the fact that customers may be routed to another queue and be served again (even multiple times) during the same cycle.

From branching theory we adopt the term \emph{descendants} of a certain (tagged) customer to denote all customers that arrive (in all queues) during the service of this tagged customer, plus the customers arriving during their service times, and so on. If, upon his service completion, a customer is routed to another queue, we also consider him as his own descendant. We define $B^*_{k,i}$, $i=1,\dots,N; k=0, \dots, N$, as the service time of a type $i-k$ (which is understood as $N+i-k$ if $i \leq k$) customer at $Q_{i-k}$, plus the service times of all of his descendants that will be served before or during the next visit of the server to $Q_i$. The special case $B^*_{0,i}$ is simply the service time of a type $i$ customer, $i=1,\dots,N$. A formal definition in terms of LSTs is given below:
\begin{align}
\Btilde_{k,i}(\omega) &= \B_{i-k}\Big(\omega + \sum_{j=0}^{k-1}\lambda_{i-j}\big(1-\Btilde_{j,i}(\omega)\big)\Big)\Ptilde_{k,i}(\omega),  &&k=0,1,\dots,N; i=1,\dots,N,
\label{bkiomega}\\
\intertext{where}
\Ptilde_{k,i}(\omega) &= 1-\sum_{j=0}^{k-1}p_{i-k,i-j}\big(1-\Btilde_{j,i}(\omega)\big), &&k=0,1,\dots,N; i=1,\dots,N.    \label{pkiomega}
\end{align}
For a type $i-k$ customer, $P^*_{k,i}$ accounts for the service times of his descendants that are caused by the fact that he may be routed to another queue upon his service completion.

A similar function should be defined for the switch-over times:
\begin{equation}
\Rtilde_{k,i}(\omega) = \R_{i-k}\Big(\omega + \sum_{j=0}^{k-1}\lambda_{i-j}\big(1-\Btilde_{j,i}(\omega)\big)\Big),  \qquad\ \qquad k=0,1,\dots,N; i=1,\dots,N.
\label{rkiomega}
\end{equation}
Note that, compared to \eqref{bkiomega}, no term $\Ptilde_{k,i}(\omega)$ is required because no routing takes place at the end of a switch-over time.

Finally, we define the following $N+1$ dimensional vectors:
\begin{align}
\ub{k,i} &= \big(1, \dots, 1, \Btilde_{k,i}(\omega), 1, \dots, 1\big), && k=0,1,\dots,N-1; i=1,\dots,N, \label{ubki}\\
\ub{N,i} &= \big(1, \dots, 1, \Btilde_{0,i}(\omega)\big), && i=1,\dots,N,\label{ubNi}
\end{align}
with $\Btilde_{k,i}(\omega)$ at position $i-k$ in \eqref{ubki} (or position $N+i-k$ if $k\geq i$), and $\Btilde_{0,i}(\omega)$ at position $N+1$ in \eqref{ubNi}. We use $\bigotimes$ to denote the element-wise multiplication of vectors.

\begin{proposition}\label{cycletimethm}
The LST of the distribution of the cycle time $C_i$ is given by
\begin{equation}
\C_i(\omega) = \LB^{(V_i)}\big( \bigotimes_{k=0}^{N-1} \ub{k,i-1}\big)\prod_{k=0}^{N-1}\Rtilde_{k,i-1}(\omega), \qquad i=1,\dots,N.\label{cycleTimeLST}
\end{equation}
The interpretation of \eqref{cycleTimeLST} is that the length of a cycle starting with a visit to $Q_i$ is the sum of the \emph{extended} service times of all customers present at the beginning of the cycle, and the sum of all \emph{extended} switch-over times during the cycle. By extended service time (switch-over time) we refer to a service time (switch-over time) plus the service times of all customers that arrive during this service time (switch-over time) in one of the queues that are yet to be served during the remainder of the cycle, and all of their descendants that will be served before the end of the cycle.
\begin{proof}
To prove Proposition \ref{cycletimethm} we keep track of all the customers that will be served during one cycle. We condition on the numbers of customers present in each queue at the beginning of $C_i$, denoted by $n_1, \dots, n_N$. Note that there are no gated customers present at this moment, because the gate has been removed at the beginning of the last switch-over time of the previous cycle. A cycle $C_i$ consists of:
\begin{enumerate}
\item the service of all customers present at the beginning of the cycle,
\item all of their descendants that will be served before the start of the next cycle (i.e., before the next visit to $Q_i$),
\item the switch-over times $R_1, \dots, R_N$,
\item all customers arriving during these switch-over times that will be served before the start of the next cycle,
\item all of their descendants that will be served before the start of the next cycle.
\end{enumerate}
We define $S_j$ for $j=1,\dots, N$, as the service time of a type $j$ customer plus the service times of all of his descendants that will be served during (the remaining part of) $C_i$. Since the service discipline is gated at all queues, we have:
\begin{equation}
S_j = B_j + \sum_{k=j+1}^{i-1}\sum_{l=1}^{N_k(B_j)} S_{k_l} +
\begin{cases}
S_m & \qquad\text{for $m=j+1,\dots,i-1$, w.p. $p_{j,m}$},\\
0 & \qquad\text{w.p. }1-\sum_{m=j+1}^{i-1}p_{j,m},
\end{cases}
\label{Sj}
\end{equation}
where $N_k(T)$ denotes the number of arrivals in $Q_k$ during a (possibly random) period of time $T$, and $S_{k_l}$ is a sequence of (independent) extended service times $S_k$. Note that $S_j$ depends on $i$, although we have chosen to hide this for presentational purposes. The gated service discipline is reflected in the fact that only customers arriving in (or rerouted to) $Q_{j+1},\dots,Q_{i-1}$ are being served during the residual part of $C_i$. It can easily be shown that the LST of $S_{i-k}$ is $\Btilde_{k-1,i-1}(\omega)$ for $k=1, \dots, N$.
Note that the first summation in \eqref{Sj} is cyclic, which may sometimes cause confusion  (for example if $j=i-1$, when this is supposed to be a summation over zero terms). Avoiding this (possible) confusion is the main reason that we have chosen to define $\Btilde_{k,i}(\omega)$, $\Ptilde_{k,i}(\omega)$ and $\Rtilde_{k,i}(\omega)$ relative to queue $i$ ($k$ steps backward in time).

Using this branching way of looking at the cycle time, we can express $C_i$ in terms of $R_1, \dots, R_N$ and $S_1, \dots, S_N$.
First, however, we derive the following intermediate result.
\begin{align*}
\E\left[\ee^{-\omega R_{i-k}}\prod_{j=i-k+1}^{i-1}\prod_{l=1}^{N_j(R_j)}\ee^{-\omega S_{j_l}}\right]&
%&=\int_{t=0}^\infty \ee^{-\omega t} \prod_{j=i-k+1}^{i-1}
%\sum_{m_j=0}^\infty\E\left[\prod_{l=1}^{m_j}\ee^{-\omega S_{j_l}}\right]\P[N_j(t_j)=m_j]\dd \P[R_{i-k} < t]\\
%&=\int_{t=0}^\infty \ee^{-\omega t_j} \prod_{j=i-k+1}^{i-1}
%\sum_{m_j=0}^\infty\E\left[\ee^{-\omega S_j}\right]^{m_j}\frac{(\lambda_j t_j)^{m_j}}{m_j!}\ee^{-\lambda_j t_j}\dd \P[R_{i-k} < t]\\
%&=\int_{t=0}^\infty  \ee^{-\left(\omega + \sum_{j=i-k+1}^{i-1}\lambda_j(1- \E[\ee^{-\omega S_j}])\right)t   }\dd \P[R_{i-k} < t]\\
=\R_{i-k}\big(\omega + \sum_{j=i-k+1}^{i-1}\lambda_j(1- \E[\ee^{-\omega S_j}])\big)\\
&=\Rtilde_{k-1,i-1}(\omega).
\end{align*}

Now, introducing the shorthand notation $n_1,\dots,n_N$ for the event that the numbers of customers at the beginning of $C_i$ in queues $1,\dots,N$ are respectively $n_1,\dots,n_N$, we can find the cycle time LST conditional on this event.

\begin{align*}
\E\left[\ee^{-\omega C_i}\,|\,n_1,\dots,n_N\right] &= \E\left[\exp\Big(-\omega
\sum_{j=i-N}^{i-1}  \big(\sum_{l=1}^{n_j} S_{j_l}         + R_j + \sum_{k=j+1}^{i-1}\sum_{l=1}^{N_k(R_j)}S_{k_l} \big)
\Big)\right]\\
&=\E\left[\prod_{j=i-N}^{i-1} \left(\prod_{l=1}^{n_j}  \ee^{-\omega S_{j_l}}\right)   \ee^{-\omega R_j}\prod_{k=j+1}^{i-1}\prod_{l=1}^{N_k(R_j)}\ee^{-\omega S_{k_l}}\right]\\
&=\prod_{j=i-N}^{i-1}\left(\prod_{l=1}^{n_j}\E\left[\ee^{-\omega S_{j_l}}\right]\right)\prod_{j=i-N}^{i-1}\E\left[\ee^{-\omega R_j}\prod_{k=j+1}^{i-1}\prod_{l=1}^{N_k(R_j)}\Big(\ee^{-\omega S_{k_l}}\Big)\right]\\
%&=\prod_{j=i-N}^{i-1}\E\left[\ee^{-\omega S_j}\right]^{n_j}\prod_{j=i-N}^{i-1}
%\int_{t_j=0}^\infty \ee^{-\omega t_j} \prod_{k=j+1}^{i-1}
%\sum_{m_k=0}^\infty\E\left[\prod_{l=1}^{m_k}\Big(\ee^{-\omega S_{k_l}}\Big)\right]\P[N_k(t_j)=m_k]\dd \P[R_j < t_j]\\
%&=\left(\prod_{k=1}^N\Btilde_{k-1, i-1}(\omega)^{n_{i-k}}\right) \prod_{j=i-N}^{i-1}
%\int_{t_j=0}^\infty \ee^{-\omega t_j} \prod_{k=j+1}^{i-1}
%\sum_{m_k=0}^\infty\E\left[\ee^{-\omega S_k}\right]^{m_k}\frac{(\lambda_k t_j)^{m_k}}{m_k!}\ee^{-\lambda_k t_j}\dd \P[R_j < t_j]\\
%&=\left(\prod_{k=1}^N\Btilde_{k-1, i-1}(\omega)^{n_{i-k}}\right) \prod_{j=i-N}^{i-1}
%\int_{t_j=0}^\infty  \ee^{-\left(\omega + \sum_{k=j+1}^{i-1}\lambda_k(1- \E[\ee^{-\omega S_k}])\right)t_j   }\dd \P[R_j < t_j]\\
%&=\left(\prod_{k=1}^N\Btilde_{k-1, i-1}(\omega)^{n_{i-k}}\right) \prod_{j=i-N}^{i-1}
%\R_j\big(\omega + \sum_{k=j+1}^{i-1}\lambda_k(1- \E[\ee^{-\omega S_k}])\big)\\
&=\left(\prod_{k=1}^N\Btilde_{k-1, i-1}(\omega)^{n_{i-k}}\right) \prod_{k=1}^N \Rtilde_{k-1,i-1}(\omega).\\
\end{align*}
Equation \eqref{cycleTimeLST} follows after deconditioning.
\end{proof}
\end{proposition}

\begin{remark}\label{alternativeC}
Because of our main interest in the waiting-time distributions, we have followed quite an elaborate path to find the LST of the cycle-time distribution. However, if one is merely interested in a quick way to find $\C_i(\omega)$, a more efficient approach can be used. One of the most efficient ways to find $\C_i(\omega)$ is to distinguish between customers that arrive from outside the network (external customers) and internally rerouted customers (internal customers). One can straightforwardly adapt the laws of motion \eqref{lawsofmotion1}-\eqref{lawsofmotionN} to find an expression for $\LB^{(V_i)'}(z_1^E, z_1^I, \dots, z_N^E, z_N^I)$.
Just like $\LB^{(V_i)}(z_1, \dots, z_N, z_G)$, $\LB^{(V_i)'}(z_1^E, z_1^I, \dots, z_N^E, z_N^I)$ stands for the PGF of the joint queue length at the beginning of $V_i$, but now we distinguish between external and internal customers in each queue (indicated by $z_j^E$ and $z_j^I$). Since external customers arrive in $Q_i$ according to a Poisson process with intensity $\lambda_i$, one can apply the distributional form of Little's Law (see, for example, Keilson and Servi \cite{keilsonservi90}) to the \emph{external} customers in $Q_i$:
\[
\C_i(\omega) = \LB^{(V_i)'}(1, \dots, 1, 1-\omega/\lambda_i, 1, \dots, 1), \qquad i=1,\dots, N.
\]

\end{remark}

\subsection{Waiting-time distributions}\label{waitingTimesubsection}

In this subsection we find the LSTs of $W_i^{E}$ and $W_i^I$, the waiting-time distributions of arbitrary external and internal customers in $Q_i$, and use them to obtain the LST of $W_i$, the waiting time of an arbitrary customer. Recall that the waiting time $W_i$ of an arbitrary customer in $Q_i$ is the time between his arrival at this queue and the moment at which his service starts. Hence, even if a customer is routed to the same queue multiple times, each visit to this queue invokes a new waiting time. We stress that common methods used in the polling literature to find waiting time LSTs cannot be applied in our queueing network, because they rely heavily on the assumption that \emph{every} customer in the system has arrived according to a Poisson process. Since this assumption is violated in our model, we have developed a novel approach to find the waiting time LST of an arbitrary customer in our network. The joint queue-length distributions at various epochs, as discussed in Subsection \ref{jointQLsubsection}, play an essential role in the analysis. First we focus on the waiting times of internal customers, then we discuss the waiting times of external customers.

\paragraph{Internal customers.}
The arrival epoch of an internal customer always coincides with a service completion. Hence, we condition on the joint queue length and the arrival epoch of an internal customer to find his waiting time LST. The waiting time of an internal customer \emph{given that} he arrives in $Q_i$ after a service completion at $Q_{i-k}$ is denoted by $\textit{WC}_{i}^{(B_{i-k})}$ ($i,k = 1,\dots,N$). To find $\textit{WC}_{i}^{(B_{i-k})}$, we only have to compute the probability that an arbitrary internal customer in $Q_i$ arrives after a service completion at $Q_{i-k}$. The mean number of customers (internal plus external) present at the beginning of $V_{i-k}$ at $Q_{i-k}$ is $\gamma_{i-k}\E[C]$. Each of these customers joins $Q_i$ upon his service completion with probability $p_{{i-k},i}$. This observation combined with the fact that the mean number of \emph{internal} customers arriving at $Q_i$ during the course of one cycle is $(\gamma_i-\lambda_i)\E[C]$, leads to the following result:
\begin{equation}
\W_i^I(\omega) = \sum_{k=1}^N \frac{\gamma_{i-k}p_{{i-k},i}}{\gamma_i-\lambda_i}\WC_{i}^{(B_{i-k})}(\omega), \qquad i=1,\dots,N.
\label{WiI}
\end{equation}
As a consequence, the problem of finding $\W_i^I(\cdot)$ is reduced to finding $\WC_{i}^{(B_{i-k})}(\omega)$ for all $i,k=1,\dots,N$.

For notational reasons we first introduce the following $N+1$ dimensional vectors, which will appear several times in this section:
\begin{align*}
\ubG{k,i} &= \begin{cases}
\displaystyle\ub{0,i}  & \qquad\textrm{ if }k<0,\\
\displaystyle\ub{0,i}\bigotimes_{j=0}^{k-1}\ub{j,i-1}  & \qquad\textrm{ if }k=1,\dots,N,\\
\displaystyle\ub{N,i}\bigotimes_{j=0}^{N-1}\ub{j,i-1}  & \qquad\textrm{ if }k=N,
\end{cases}
\end{align*}
for $i=1,\dots,N$. Again, we use $\bigotimes$ to denote the element-wise multiplication of vectors.

\begin{proposition}\label{LSTWexternalthm}
We have
\begin{align}
\WC_{i}^{(B_{i-k})}(\omega) &= \LC^{(B_{i-k})}\big( \ubG{k,i}\big)\prod_{j=0}^{k-1}\Rtilde_{j,i-1}(\omega),
\label{WkiI}
\end{align}
for $i,k=1,\dots,N$.
\begin{proof}
The key observation in the proof of Proposition \ref{LSTWexternalthm} is that an arrival of an internally rerouted customer always coincides with some service completion. For this reason, we consider the system right after the service completion at, say, $Q_j$ ($j=1,\dots,N$). We compute the waiting time LST of a customer routed to $Q_i$ after being served in $Q_j$, conditional on the numbers of customers of each type (now \emph{including} gated customers) present at the arrival epoch (\emph{not} including the arriving customer himself). We denote by $n_1, \dots, n_N, n_G$ the event that the numbers of customers of all types are respectively $n_1, \dots, n_N, n_G$. Let $n_{iG}:=n_i$ if $i\neq j$, and $n_{iG}:=n_G$ if $i= j$. Note that the type $G$ customers are located behind the gate in $Q_j$, and that the customer routed to $Q_i$ only has to wait for these customers in case $i=j$. The waiting time of the tagged customer consists of:
\begin{enumerate}
\item the service of all $n_j$ customers in front of the gate in $Q_j$ at the arrival epoch,\label{enum1}
\item the service of all $n_{j+1},\dots,n_{i-1}$ customers present in $Q_{j+1},\dots,Q_{i-1}$ at the arrival epoch,\label{enum2}
\item all of the descendants of the previously mentioned customers that will be served before the next visit to $Q_i$,
\item if $i\neq j$, the service of all $n_{iG}$ customers present in $Q_i$ at the arrival epoch; if $i=j$, the service of all $n_{iG}$ gated customers present in $Q_i$ at the arrival epoch,
\item the switch-over times $R_j, \dots, R_{i-1}$,
\item all customers arriving during these switch-over times that will be served before the next visit to $Q_i$,
\item all of their descendants that will be served before the next visit to $Q_i$.
\end{enumerate}
We denote the  waiting time of an internal customer conditional on the event that he arrives in $Q_i$ after being served in $Q_j$, \emph{and conditional on the event that the numbers of customers of all types at the arrival epoch are respectively} $n_1, \dots, n_N, n_G$, by $\textit{WC}_i^{(B_j)'}$. Just like in the proof of Proposition \ref{cycletimethm}, we can express $\textit{WC}_i^{(B_j)'}$ in terms of $R_1, \dots, R_N$ and $S_1, \dots, S_N$:
\begin{align}
\textit{WC}_i^{(B_j)'} &= \sum_{k=j}^{i-1}\left[\sum_{l=1}^{n_k}S_{k_l} + R_k+\sum_{l=k+1}^{i-1}\sum_{m=1}^{N_l(R_k)}S_{l_m}\right] + \sum_{l=1}^{n_{iG}}B_{i,l}.\label{WCiBj}
\end{align}
Taking the LST of \eqref{WCiBj} leads to \eqref{WkiI} after deconditioning. The derivation proceeds along the exact same lines as in the proof of Proposition \ref{cycletimethm}, and is therefore omitted.

\end{proof}
\end{proposition}

\paragraph{External customers.}
External customers arrive in $Q_i$ according to a Poisson process with intensity $\lambda_i$. We distinguish between customers arriving during a switch-over time and customers arriving during a visit time. The waiting time of an external customer in $Q_i$ \emph{given that} he arrives during $R_{i-k}$ is denoted by $W_{i}^{(R_{i-k})}$ ($i,k = 1,\dots,N$). Similarly, we use $W_{i}^{(V_{i-k})}$ to denote an external customer arriving in $Q_i$ during $V_{i-k}$. The waiting time LST of an arbitrary external customer can be expressed in terms of $\W_{i}^{(R_{i-k})}(\cdot)$ and $\W_{i}^{(V_{i-k})}(\cdot)$:
\begin{equation}
\W_i^E(\omega) = \frac{1}{\E[C]}\sum_{k=1}^N \left(\E[V_{i-k}]\W_{i}^{(V_{i-k})}(\omega)+r_{i-k}\W_{i}^{(R_{i-k})}(\omega)\right), \qquad i=1,\dots,N.
\label{WiE}
\end{equation}
We  first focus on the waiting time of customers arriving during a switch-over time. Consider a tagged customer arriving in $Q_i$ during $R_{i-k}$, $i,k=1,\dots,N$. Since the remaining part of the switch-over time is part of the waiting time of the arriving customer, it will turn out that we need the \emph{joint} distribution of all customers present at the arrival epoch \emph{and} the residual part of $R_{i-k}$, denoted by $R_{i-k}^R$. The PGF of the joint queue-length distribution at the arrival epoch is given by \eqref{jointQLduringSj}. Equation \eqref{jointQLduringSj} is based on the observation that the number of customers in each queue at an arbitrary moment during $R_{i-k}$ is simply the sum of the number of customers present at the beginning of $R_{i-k}$ and the number of customers that have arrived during the elapsed (past) part of $R_{i-k}$, denoted by $R_{i-k}^P$. These random variables are independent. Hence, it is straightforward to adapt \eqref{jointQLduringSj} to find the joint distribution of the queue lengths \emph{and} residual part of $R_{i-k}$, using the following result from elementary renewal theory:
\[
\R^{PR}_j(\omega_P, \omega_R) = \frac{\R_j(\omega_P)-\R_j(\omega_R)}{(\omega_R-\omega_P)r_j},\qquad j=1,\dots,N,
\]
with $\R^{PR}_j(\omega_P, \omega_R)$ denoting the LST of the joint distribution of past and residual switch-over time $R_j$. Hence,
\begin{equation}
%\E\big[z_1^{\textit{LB}_1^{(R_j)}} \dots z_N^{\textit{LB}^{(R_j)}_N} z_G^{\textit{LB}_G^{(R_j)}} \ee^{-\omega_R R_j^R}\big]
\L^{(R_j)}(\z, \omega)= \LB^{(R_j)}(\z)\R^{PR}_j(\Sigma(\z), \omega),
\label{jointQLduringSjandRjres}
\end{equation}
where $\L^{(R_j)}(\z, \omega)$ denotes the PGF-LST of the joint distribution of the number of customers of each type at an arbitrary moment during $R_j$ and the residual part of $R_j$. Obviously, there are no gated customers present during a switch-over time.

Consequently, and also using PASTA, we can find the waiting-time distribution by conditioning on the number of customers present at an arbitrary moment during $R_{i-k}$ and on the residual switch-over time.

\begin{proposition}
We have
\begin{align}
\W_{i}^{(R_{i-k})}(\omega) &= \R^{PR}_{i-k}\Big(\sum_{j=1}^{k-1}\lambda_{i-j}\big(1-\Btilde_{j-1,i-1}(\omega)\big) + \lambda_{i}\big(1-\B_i(\omega)\big), \omega+\sum_{j=1}^{k-1}\lambda_{i-j}\big(1-\Btilde_{j-1,i-1}(\omega)\big)\Big)\nonumber\\
&\times \LB^{(R_{i-k})}\big( \ubG{k-1,i}\big)
\prod_{j=0}^{k-2}\Rtilde_{j,i-1}(\omega),  \qquad i,k=1,\dots,N,
\label{WkiR}
\end{align}
\begin{proof}
We consider an arbitrary customer arriving in $Q_i$ during $R_j$. Similar to the proofs of the preceding propositions in this section, we condition on the number of customers present in all queues at the arrival epoch, denoted by $n_1, \dots, n_N$. As mentioned before, no gated customers are present during a switch-over time. However, we also condition on the residual length of $R_j$, denoted by $t_R$. The waiting time of the tagged customer consists of:
\begin{enumerate}
\item the service of all $n_{j+1},\dots,n_{i-1}$ customers present at the arrival epoch in $Q_{j+1},\dots,Q_{i-1}$,
\item the service of all their descendants that will be served before the start of the next visit to $Q_i$,
\item the service of all $n_i$ customers present at the arrival epoch in $Q_i$,
\item the residual switch-over time $t_R$,
\item the switch-over times $R_{j+1}, \dots, R_{i-1}$,
\item the service of all customers arriving during $t_R, R_{j+1}, \dots, R_{i-1}$ that will be served before the start of the next visit to $Q_i$,
\item the service of all descendants of these customers that will be served before the start of the next visit to $Q_i$.
\end{enumerate}
If we denote the waiting time of a type $i$ customer arriving during $R_j$, \emph{conditional on} $n_1, \dots, n_N$ and $t_R$, by $W_i^{(R_j)'}$, we can summarise these items in the following formula:
\begin{align}
W_i^{(R_j)'} &= \sum_{k=j+1}^{i-1}\left[\sum_{l=1}^{n_k}S_{k_l} + R_k + \sum_{l=k+1}^{i-1}\sum_{m=1}^{N_l(R_k)}S_{l_m}\right]+\sum_{l=1}^{n_i}B_{i_l}
+t_R + \sum_{l=j+1}^{i-1}\sum_{m=1}^{N_l(t_R)}S_{l_m}.\label{WiRj}
\end{align}
Taking the LST of \eqref{WiRj} and using \eqref{jointQLduringSjandRjres} leads to \eqref{WkiR} after deconditioning. The derivation is not completely straightforward, but rather than providing it here, we refer to the proof of Proposition \ref{WkiVthm}, which contains a similar derivation of a more complicated equation.
\end{proof}
\end{proposition}

Now we only need to determine $\W_{i}^{(V_{i-k})}(\cdot)$. Focussing on a tagged customer arriving in $Q_i$ during the service of a customer in $Q_{i-k}$, for $i,k=1,\dots,N$, we can find $\W_{i}^{(V_{i-k})}(\cdot)$ by conditioning on the number of customers in each queue at the arrival epoch and the residual service time. Similar to $\R^{PR}_j(\cdot)$, we define the LST of the joint distribution of past and residual service time $B_j$ as
\begin{equation}
\B^{PR}_j(\omega_P, \omega_R) = \frac{\B_j(\omega_P)-\B_j(\omega_R)}{(\omega_R-\omega_P)b_j},\qquad j=1,\dots,N.
\label{BPRj}
\end{equation}
%Similar to what we have done before for $\L^{(R_j)}(\z, \omega)$,
We can now use Equations \eqref{jointQLduringVj} and \eqref{BPRj} to find the PGF-LST of the joint distribution of the number of customers of each type present at an arbitrary moment during $V_j$ and the residual service time of the customer that is being served at that moment:
\begin{equation}
\L^{(V_j)}(\z, \omega)= \LB^{(B_j)}(\z)\B^{PR}_j(\Sigma_j(\z), \omega).
\end{equation}
Note that the customers arriving in $Q_j$ during the elapsed part of $B_j$ are gated customers.

\begin{proposition}\label{WkiVthm}
We have
\begin{align}
\W_{i}^{(V_{i-k})}(\omega) &=
\B^{PR}_{i-k}\Big(\sum_{j=1}^{k-1}\lambda_{i-j}\big(1-\Btilde_{j-1,i-1}(\omega)\big) + \lambda_{i}\big(1-\B_i(\omega)\big), \omega+\sum_{j=1}^{k-1}\lambda_{i-j}\big(1-\Btilde_{j-1,i-1}(\omega)\big)\Big)\nonumber\\
&\times \LB^{(B_{i-k})}\big( \ubG{k,i}\big)
\prod_{j=0}^{k-1}\Rtilde_{j,i-1}(\omega)%\nonumber\\
\times\frac{\Ptilde_{k-1,i-1}(\omega)}{\Btilde_{k-1,i-1}(\omega)},
\label{WkiV1}
\end{align}
for $i,k=1,\dots,N$.
\begin{proof}
We denote by $n_1, \dots, n_N, n_G$ the numbers of customers of all types present at the arrival epoch of the tagged customer. The residual part of the service time of the customer being served at this arrival epoch is denoted by $t_R$. Let $n_{iG}:=n_i$ if $i\neq j$, and $n_{iG}:=n_G$ if $i= j$. The waiting time of a type $i$ customer arriving during $V_j$, conditional on $n_1, \dots, n_N, n_G$ and the residual service time consists of the following components:
\begin{enumerate}
\item the service of $n_j-1$ customers in front of the gate in $Q_j$ (We exclude the customer being served at the arrival epoch),
\item the service of all $n_{j+1},\dots,n_{i-1}$ customers present in $Q_{j+1},\dots,Q_{i-1}$,
\item all of the descendants of the previously mentioned customers that will be served before the next visit to $Q_i$,
\item if $i\neq j$, the service of all $n_{iG}$ customers present in $Q_i$ at the arrival epoch; if $i=j$, the service of all $n_{iG}$ gated customers present in $Q_i$,
\item the switch-over times $R_j, \dots, R_{i-1}$,
\item the residual service time $t_R$,
\item all customers arriving during $t_R$ and $R_j, \dots, R_{i-1}$ that will be served before the next visit to $Q_i$,
\item all of their descendants that will be served before the next visit to $Q_i$,
\item the (possible) future service of the customer being served at the arrival epoch, due to the fact that he may be routed to another queue that will be served before the next visit to $Q_i$,
\item the service of all descendants of this rerouted customer (Note that if he will be rerouted and served again, he will count as his own descendant).
\end{enumerate}
More formally:
\begin{equation}
\begin{aligned}
W_i^{(V_j)'} &= \sum_{l=1}^{n_j-1}S_{j,l} + \sum_{k=j+1}^{i-1}\sum_{l=1}^{n_k}S_{k_l}+\sum_{l=1}^{n_{iG}}B_{i_l} +\sum_{k=j}^{i-1}\left[R_k + \sum_{l=k+1}^{i-1}\sum_{m=1}^{N_l(R_k)}S_{l_m}\right]  \\
&+t_R + \sum_{l=j+1}^{i-1}\sum_{m=1}^{N_l(t_R)}S_{l_m}+\begin{cases}
S_l & \qquad\text{for $l=j+1,\dots,i-1$, w.p. $p_{j,l}$},\\
0 & \qquad\text{w.p. }1-\sum_{l=j+1}^{i-1}p_{j,l},
\end{cases}.
\end{aligned}
\label{WiVj}
\end{equation}
We now show that Equation \eqref{WkiV1} follows from taking the LST:
\begin{align*}
&\E[\ee^{-\omega W_i^{(V_j)}}|n_1,\dots,n_N,n_{iG}] \\
&=\E\left[\prod_{l=1}^{n_j-1}\ee^{-\omega S_{j_l}}\prod_{m=j+1}^{i-1}\prod_{l=1}^{n_m}\ee^{-\omega S_{m_l}}\right]\E\left[\prod_{l=1}^{n_{iG}}\ee^{-\omega B_{i_l}}\right]\E\left[\prod_{m=j}^{i-1}\ee^{-\omega\left(R_m + \sum_{l=m+1}^{i-1}\sum_{q=1}^{N_l(R_m)}S_{l_q}  \right)}\right]\\
&\times \ee^{-\omega t_R} \E\left[\prod_{l=j+1}^{i-1}\prod_{m=1}^{N_l(t_R)}\ee^{-\omega S_{l_m}}    \right]
\left(\sum_{l=j+1}^{i-1}p_{j,l}\E\left[\ee^{-\omega S_l}\right]+1-\sum_{l=j+1}^{i-1}p_{j,l}\right)\\
&= \E\left[\ee^{-\omega S_{j}}\right]^{n_j-1}\prod_{m=j+1}^{i-1}\E\left[\ee^{-\omega S_{m}}\right]^{n_m}\E\left[\ee^{-\omega B_{i}}\right]^{n_{iG}}
\prod_{m=j}^{i-1}\R_m\Big( \omega + \sum_{l=m+1}^{i-1}(1-\E[\ee^{-\omega S_l}])\Big)\\
&\times \ee^{-\omega t_R} \prod_{l=j+1}^{i-1}\sum_{m=0}^{\infty}\E[\ee^{-\omega S_{l}}]^{m}\P[N_l(t_R)=m]
\left(1-\sum_{l=j+1}^{i-1}p_{j,l}\Big(1-\E\left[\ee^{-\omega S_l}\right]\Big)\right)\\
&= \Btilde_{k-1,i-1}(\omega)^{n_{i-k}-1}\prod_{l=1}^{k-1}\Btilde_{l-1,i-1}(\omega)^{n_{i-l}}\B_{i}(\omega)^{n_{iG}}
\prod_{l=1}^{k}\Rtilde_{l-1,i-1}(\omega)\\
&\times \exp\left[-\Big(\omega + \sum_{l=j+1}^{i-1}(1-\E[\ee^{-\omega S_{l}}])\Big)t_R\right] \Ptilde_{k-1,i-1}(\omega)\\
&= \Btilde_{k-1,i-1}(\omega)^{n_{i-k}}\prod_{l=1}^{k-1}\Btilde_{l-1,i-1}(\omega)^{n_{i-l}}\B_{i}(\omega)^{n_{iG}}
\prod_{l=1}^{k}\Rtilde_{l-1,i-1}(\omega)\\
&\times \exp\left[-\Big(\omega + \sum_{l=1}^{k-1}(1-\Btilde_{l-1,i-1}(\omega))\Big)t_R\right] \frac{\P_{k-1,i-1}(\omega)}{\Btilde_{k-1,i-1}(\omega)},\\
\end{align*}
where $k=i-j$ (or $k=N+i-j$ if $j\geq i$).
Deconditioning of this expression leads to \eqref{WkiV1}.
\end{proof}
\end{proposition}

\paragraph{Arbitrary customers.}
Finally, we present the main result of this section: the LST of the waiting-time distribution of an arbitrary customer in $Q_i$.
\begin{theorem}
The LST of the waiting-time distribution of an arbitrary customer in $Q_i$, if this queue receives gated service, is given by:
\begin{equation}
\W_i(\omega) = \frac{\gamma_i-\lambda_i}{\gamma_i}\W_i^I(\omega) + \frac{\lambda_i}{\gamma_i}\W_i^E(\omega), \qquad i=1,\dots,N,\label{LSTW}
\end{equation}
where $\W_i^I(\omega)$ and $\W_i^E(\omega)$ are given by \eqref{WiI} and \eqref{WiE}, respectively.
\begin{proof}
The result follows immediately after conditioning on the event that an arbitrary customer is an internal or external customer.
\end{proof}
\end{theorem}

\section{Exhaustive service}\label{exhaustiveservice}

In this section we study systems with mixtures of gated and exhaustive service, that is, some queues are served exhaustively whereas other queues receive gated service. We restrict ourselves to presenting the results, but for reasons of compactness we omit all proofs as they can be produced similar to the proofs in the previous section.

Throughout we use the index $e\in\{1,\dots,N\}$ to refer to an arbitrary queue with exhaustive service, which means that customers are being served until the queue is empty. This means that, in contrast to gated service, customers arriving in $Q_e$ \emph{during} $V_e$ will be served during that same visit period. This is true, even if the customer has just received service in $Q_e$ and was routed back to $Q_e$ again. To deal with this issue, we define an extended service time $B_e^\textit{exh}$ which is the total amount of service that a customer receives during a visit period $V_e$ before being routed to another queue (or leaving the system), cf. \cite{sidi2}. As stated in \cite{sidi2}, $B_e^\textit{exh}$ is the geometric sum, with parameter $p_{e,e}$, of independent random variables with the same distribution as $B_e$. The LST of $B_e^\textit{exh}$ is given by
\[
\Bexh_e(\omega) = \frac{(1-p_{e,e})\B_e(\omega)}{1-p_{e,e}\B_e(\omega)}.
\]
We denote a busy period of type $e$ customers by $\textit{BP}_e$. The PGF-LST of the joint distribution of a busy period and the number of customers served during this busy period satisfies the following equation:
\[
\BP_e(z, \omega) = z\Bexh_e\big(\omega + \lambda_e(1-\BP_e(z, \omega))\big).
\]

\subsection{The joint queue-length distributions}
\paragraph{Visit beginnings and completions.}
The laws of motion \eqref{lawsofmotion1}-\eqref{lawsofmotionN} have to be adapted if a queue receives exhaustive service. First we need to redefine $\Sigma_i(\z)$ and $P_i(\z)$ if $Q_i$ is served exhaustively, and introduce $\Pexh_i(\z)$:
\begin{align*}
\Sigma_e(\z)&=\sum_{j\neq e}\lambda_j(1-z_j),\\
P_e(\z) &= p_{e,0}+\sum_{j=1}^Np_{e,j}z_j,\\
\Pexh_e(\z) &= \frac{p_{e,0}}{1-p_{e,e}}+\sum_{j\neq e}\frac{p_{e,j}}{1-p_{e,e}}z_j,
\end{align*}
for all $e\in\{1,\dots,N\}$ corresponding to queues with exhaustive service. The laws of motion now change accordingly:
\begin{align*}
\LC^{(V_e)}(\z) &= \LB^{(V_e)}\Big(z_1, \dots, z_{e-1},\BP_e\big(\Pexh_e(\z), \Sigma_e(\z)\big), z_{e+1}, \dots, z_N, 1\Big),\\
\LB^{(R_e)}(\z) &= \LC^{(V_e)}(\z),
\end{align*}
for any exhaustively served $Q_e$.

\paragraph{Service beginnings and completions.}

Eisenberg's relation \eqref{eisenberg} remains valid for queues with exhaustive service. Note that $P_e(\z)$ should \emph{not} be replaced by $\Pexh_e(\z)$ for exhaustive queues in \eqref{eisenberg}! Relation \eqref{servicebeginningscompletions} should be slightly changed for queues with exhaustive service, since customers are not placed behind a gate:
\begin{equation*}
\LC^{(B_e)}(\z) = \LB^{(B_e)}(\z)\B_e\big(\Sigma(\z)\big)/z_e.  %\label{servicebeginningscompletions}
\end{equation*}

\paragraph{Arbitrary moments.}
Equation \eqref{Lz} for the PGF of the joint queue-length distribution at arbitrary moments remains valid if some of the queues have exhaustive service. However,
$\L^{(V_j)}(\z)$ should be adapted for queues with exhaustive service by replacing gated customers with ``ordinary'' type $e$ customers:
\[
\L^{(V_e)}(\z) = \LB^{(B_e)}(\z)\frac{1-\B_e\big(\Sigma(\z)\big)}{b_e\Sigma(\z)}.
\]

\subsection{Cycle-time distributions}

The fact that customers arriving in an exhaustively served queue, say $Q_{i-k}$, during $V_{i-k}$ are served before the end of this visit period, requires changes in the definition of $\Btilde_{k,i}(\omega)$.
\begin{align}
\Btilde_{k,i}(\omega) &= \BP_{i-k}\Big(\Ptilde_{k,i}(\omega), \omega + \sum_{j=0}^{k-1}\lambda_{i-j}\big(1-\Btilde_{j,i}(\omega)\big)\Big),  &&k=0,1,\dots,N; i=1,\dots,N,\\
\intertext{where}
\Ptilde_{k,i}(\omega) &= 1-\sum_{j=0}^{k-1}\frac{p_{i-k,i-j}}{1-p_{i-k,i-k}}\big(1-\Btilde_{j,i}(\omega)\big), &&k=0,1,\dots,N; i=1,\dots,N.
\end{align}
Given this modified definition of $\Btilde_{k,i}(\omega)$, the function $\Rtilde_{k,i}(\omega)$ remains unchanged. The expression for the LST of the cycle time $C_i$, given by \eqref{cycleTimeLST}, also remains valid for systems containing exhaustively served queues.

\subsection{Waiting-time distributions}

\paragraph{Internal customers.}
The waiting time LST of internal customers \eqref{WiI} is determined by conditioning on the event that an arrival in $Q_i$ follows a service completion in some $Q_{i-k}$. As stated before, for queues with exhaustive service we need to take into account that customers that are routed back to the same queue will be served during the same visit period. For an arbitrary exhaustively served queue $Q_e$, this results in
\begin{equation}
\W_e^I(\omega) = \sum_{k=0}^{N-1} \frac{\gamma_{e-k}p_{{e-k},e}}{\gamma_e-\lambda_e}\WC_{e}^{(B_{e-k})}(\omega).\label{WeI}
\end{equation}
Compared to \eqref{WiI}, the summation starts at $k=0$ and runs up to $k=N-1$. We now introduce
\[
\ubexh{0,i}=\big(1, \dots, 1, \B_{i}(\omega), 1, \dots, 1\big),\qquad i=1,\dots,N,
\]
with $\B_{i}(\omega)$ at the position corresponding to customers in $Q_i$. If $Q_i$ has exhaustive service, there is a subtle difference with $\ub{0,i}$ which has $\BP_{i}(1,\omega)$ at position $i$. We can now determine $\WC_{e}^{(B_{e-k})}(\omega)$ for any $Q_e$ that receives exhaustive service:
\begin{align*}
\WC_{e}^{(B_{e-k})}(\omega) &= \LC^{(B_{e-k})}\big( \ubexh{0,e}\bigotimes_{j=0}^{k-1}\ub{j,e-1} \big)\prod_{j=0}^{k-1}\Rtilde_{j,e-1}(\omega),   \qquad k=1,\dots,N-1,
\\
\WC_{e}^{(B_{e})}(\omega) &= \LC^{(B_{e})}\big( \ubexh{0,e}\big).
\end{align*}
For each $Q_i$ that receives gated service, we can still use \eqref{WiI} with the modified definition of $\Btilde_{k,i}(\omega)$ for each $Q_{i-k}$ which receives exhaustive service.

\paragraph{External customers.}
The waiting time LST of external customers \eqref{WiE} is determined by conditioning on the event that an arrival in $Q_i$ takes place during $V_{i-1},\dots,V_{i-N}$ or during $R_{i-1},\dots,R_{i-N}$. Before discussing the waiting times of external customers arriving in an exhaustively served queue, it is important to realise that allowing some queues to have exhaustive service will now also require some changes to waiting times of customers arriving in a queue with gated service. This means that \eqref{WkiV1} should now become
\begin{align}
\W_{i}^{(V_{i-k})}(\omega) &=
\B^{PR}_{i-k}\Big(\sum_{j=1}^{k-1}\lambda_{i-j}\big(1-\Btilde_{j-1,i-1}(\omega)\big) + \lambda_{i}\big(1-\B_i(\omega)\big)+\lambda_{i-k}\big(1-\Btilde_{k-1,i-1}(\omega)\big), \nonumber\\ &\qquad\qquad\omega+\sum_{j=1}^{k-1}\lambda_{i-j}\big(1-\Btilde_{j-1,i-1}(\omega)\big)+\lambda_{i-k}\big(1-\Btilde_{k-1,i-1}(\omega)\big)\Big)\nonumber\\
&\times \LB^{(B_{i-k})}\big( \ub{0,i}\bigotimes_{j=0}^{k-1}\ub{j,i-1}  \big)
\prod_{j=0}^{k-1}\Rtilde_{j,i-1}(\omega)%\nonumber\\
\times\frac{1-\sum_{j=0}^{k-1}p_{i-k,i-j-1}\big(1-\Btilde_{j,i-1}(\omega)\big)}{\Btilde_{k-1,i-1}(\omega)},
\label{WkiV1exh}
\end{align}
if $Q_{i-k}$ receives exhaustive service (and $Q_i$ receives gated service). Compared to \eqref{WkiV1} we can see that there are two additional terms $\lambda_{i-k}\big(1-\Btilde_{k-1,i-1}(\omega)\big)$ which take into account that customers arriving in $Q_{i-k}$ during the elapsed \emph{and} during the residual part of the present service time $B_{i-k}$ will be served during the present visit period. Furthermore, we can see that $\Ptilde_{k-1,i-1}(\omega)$ has been replaced by $1-\sum_{j=0}^{k-1}p_{i-k,i-j-1}\big(1-\Btilde_{j,i-1}(\omega)\big)$, which is required because the customer being served should be allowed to return to $Q_{i-k}$ upon his service completion.

If $Q_e$ receives exhaustive service we have to make some additional changes. We have
\begin{equation}
\W_e^E(\omega) = \frac{1}{\E[C]}\sum_{k=1}^N \left(\E[V_{e-k+1}]\W_{e}^{(V_{e-k+1})}(\omega)+r_{e-k}\W_{e}^{(R_{e-k})}(\omega)\right),\label{WeE}
\end{equation}
where we have chosen to denote the waiting time LST of customers arriving in $Q_e$ during $V_e$ as $\W_{e}^{(V_{e})}(\omega)$ rather than $\W_{e}^{(V_{e-N})}(\omega)$ to illustrate the fact that they will be served during the same visit period. The expression for $\W_{e}^{(R_{e-k})}(\omega)$, given by \eqref{WkiR}, should be slightly modified if $Q_e$ receives exhaustive service. However, since the only required modification is that $\ub{0,i}$ should be replaced by $\ubexh{0,i}$, we refrain from giving the complete expression.

If $k>0$, the expression for $\W_{e}^{(V_{e-k})}(\omega)$ remains almost the same as \eqref{WkiV1} if $Q_{e-k}$ receives gated service, or \eqref{WkiV1exh} if $Q_{e-k}$ receives exhaustive service. The only change is, once again, that $\ub{0,i}$ should be replaced by $\ubexh{0,i}$. The case $k=0$ results in a much simpler expression, since we only have to wait for the service times of the customers that were present at the beginning of the present service (excluding the customer in service) plus the service times of the customers that have arrived in $Q_e$ during the elapsed part of the present service, plus the residual service time:
\[
\W_{e}^{(V_{e})}(\omega) = \B^{PR}_{e}\Big(\lambda_{e}\big(1-\B_e(\omega)\big), \omega\Big)\frac{\LB^{(B_{e})}\big( \ubexh{0,e}\big)}{\B_e(\omega)}.
\]

\paragraph{Arbitrary customers.}

The LST of the waiting-time distribution of an arbitrary customer in an exhaustively served queue immediately follows after conditioning on the event that an arbitrary customer is either an internal or an external customer, similar to the derivation of~\eqref{LSTW}. The result is presented in the theorem below.

\begin{theorem}
The LST of the waiting-time distribution of an arbitrary customer in $Q_i$, if this queue receives exhaustive service, is given by:
\begin{equation}
\W_i(\omega) = \frac{\gamma_i-\lambda_i}{\gamma_i}\W_i^I(\omega) + \frac{\lambda_i}{\gamma_i}\W_i^E(\omega), \qquad i=1,\dots,N,\label{LSTWexhaustive}
\end{equation}
where $\W_i^I(\omega)$ and $\W_i^E(\omega)$ are defined in \eqref{WeI} and \eqref{WeE}.

\end{theorem}

\section{Applicability of the model}

In this section we give some numerical examples that indicate the versatility of the model that we have discussed. To this end, we use some examples that can be found in the existing literature, and show how our model can be used to describe the various systems and find the relevant performance measures. Hence, most of the results presented in this section are not novel, but the way of deriving them is new.

\paragraph{Example 1: tandem queues with parallel queues in the first stage.}

\begin{figure}[ht]
\begin{center}
\includegraphics[width=0.7\linewidth]{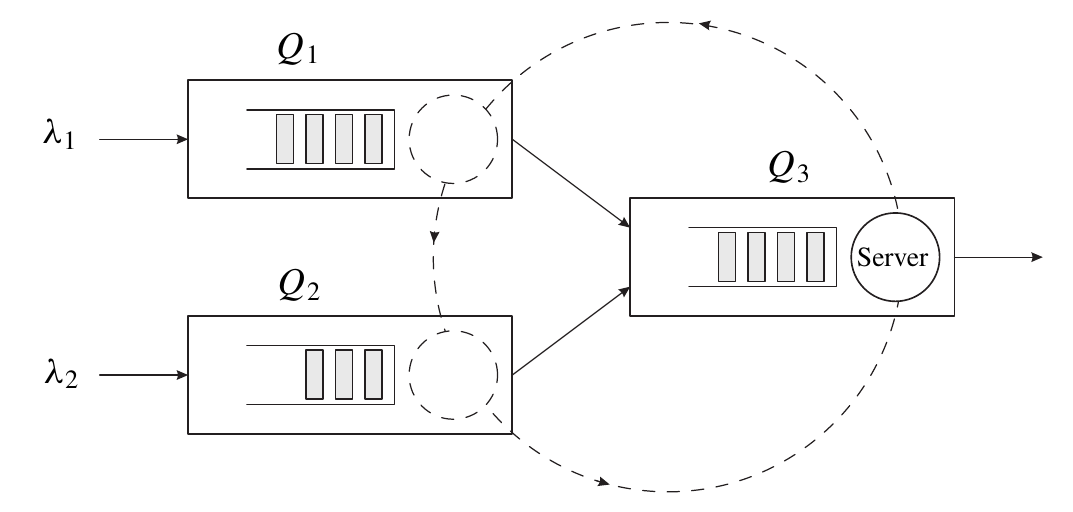}
\end{center}
\caption{Tandem queues with parallel queues in the first stage, as discussed in Example 1.}
\label{figureExample1}
\end{figure}

\begin{figure}[h]
\includegraphics[width=0.49\textwidth]{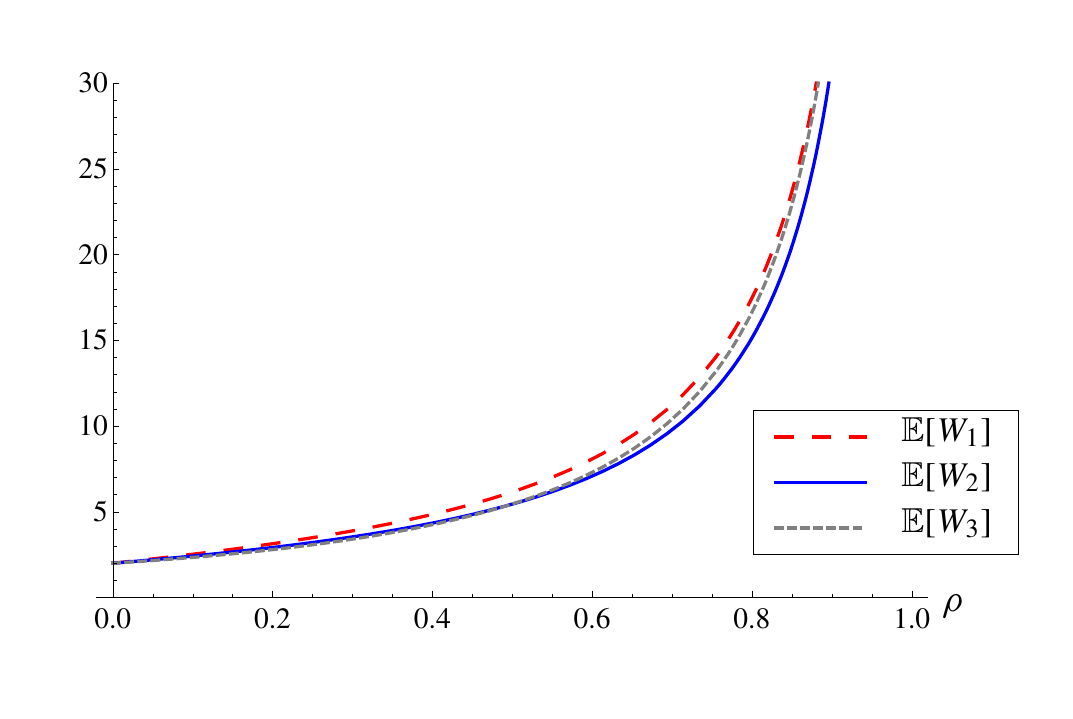}
\hfill
\includegraphics[width=0.49\textwidth]{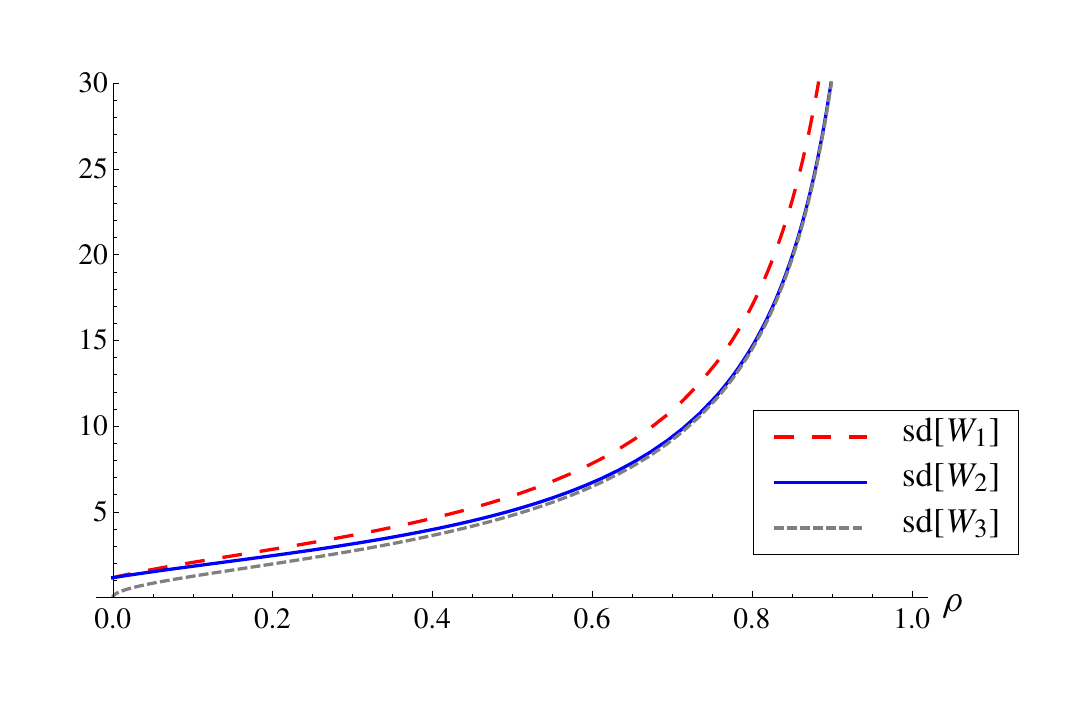}
\caption{Means and standard deviations of the waiting times in the first numerical example.}
\label{numericalresults}
\end{figure}

We first use an example that was introduced by Katayama \cite{katayama}, who studies a network consisting of three queues. Customers arrive at $Q_1$ and $Q_2$, and are routed to $Q_3$ after being served (see Figure \ref{figureExample1}). This model, which is referred to as a tandem queueing model with parallel queues in the first stage, is a special case of the model discussed in the present paper. We simply put $p_{1,3}=p_{2,3}=p_{3,0}=1$ and all other $p_{i,j}$ are zero. We use the same values as in \cite{katayama}: $\lambda_1=\lambda_2/10$, service times are deterministic with $b_1=b_2=1$, and $b_3=5$. The server serves the queues exhaustively, in cyclic order: 1, 2, 3, 1, \dots. The only difference with the model discussed in \cite{katayama} is that we introduce (deterministic) switch-over times $R_2=R_3=2$. We assume that no time is required to switch between the two queues in the first stage, so $r_1=0$. In Figure~\ref{numericalresults} we show the means and standard deviations of the waiting times of customers at the three queues. These plots reveal that in the heavy-traffic regime, as $\rho\uparrow1$, the \emph{mean} waiting times of customers in $Q_3$ are close to those in $Q_1$, but the \emph{standard deviations} of the waiting times in $Q_3$ are closer to those in $Q_2$. Further inspection of the exact results, obtained by differentiating the LSTs, confirms that in both cases the limits are very close, but not exactly the same.

It is also interesting to study the light-traffic behaviour of the system, i.e., as $\rho\downarrow0$. From the plots in Figure~\ref{numericalresults} we can see that, as $\rho\downarrow0$, the \emph{mean} waiting times are all equal, but the \emph{standard deviation} of the waiting times in $Q_1$ and $Q_2$ is different than in $Q_3$. From the LSTs of the waiting-time distributions we can obtain the exact expressions when $\rho\downarrow0$, by taking the Taylor expansion in $\rho$ at $\rho=0$ and subsequently ignoring all $\O(\rho)$ terms. This, combined with the fact that $R_1=0$ and all of the routing probabilities are either 0 or 1,  considerably simplifies all expressions from the previous section:
\begin{align*}
\W_1(\omega) &= \W_1^E(\omega) \rightarrow \frac{r_2}{r}\W_{1}^{(R_{2})}(\omega)+\frac{r_3}{r}\W_{1}^{(R_{3})},\\
\W_2(\omega) &= \W_2^E(\omega) \rightarrow \frac{r_2}{r}\W_{2}^{(R_{2})}(\omega)+\frac{r_3}{r}\W_{2}^{(R_{3})},\\
\W_3(\omega) &= \W_3^I(\omega) \rightarrow \frac{\lambda_1}{\lambda_1+\lambda_2}\WC_{3}^{(B_{1})}(\omega)+\frac{\lambda_2}{\lambda_1+\lambda_2}\WC_{3}^{(B_{2})}(\omega).
\end{align*}
Since we are considering the case $\rho\downarrow0$, these expressions can be simplified even further to closed-form expressions, because ignoring all $\O(\rho)$ terms is equivalent to regarding the system as being empty all the time:
\begin{align*}
\W_1(\omega) &\rightarrow \frac{r_2}{r}\R^{PR}_2(0, \omega)\R_3(\omega)+\frac{r_3}{r}\R^{PR}_3(0, \omega),\\
\W_2(\omega) &\rightarrow \frac{r_2}{r}\R^{PR}_2(0, \omega)\R_3(\omega)+\frac{r_3}{r}\R^{PR}_3(0, \omega),\\
\W_3(\omega) &\rightarrow \frac{\lambda_1}{\lambda_1+\lambda_2}\R_1(\omega)\R_2(\omega)+\frac{\lambda_2}{\lambda_1+\lambda_2}\R_2(\omega).
\end{align*}
These expressions reveal the true behaviour of the system in light traffic. The waiting times in $Q_1$ and $Q_2$ are simply the total residual switch-over time, with mean $r^{(2)}/2r = 2$ and second moment $r^{(3)}/3r = 16/3$. For queue $Q_3$ the situation is different, because this queue only contains internally rerouted customers. Customers being rerouted from $Q_1$ have to wait for the switch-over times $R_1+R_2$, whereas customers arriving from $Q_2$ have to wait only for $R_2$. Since $R_1=0$, the waiting time only consists of $R_2=2$ in both cases. Substituting all parameter values results in the following LT limits of the waiting-time LSTs:
\[
\W_1(\omega) \rightarrow \frac{1-\ee^{-4\omega}}{4\omega}, \qquad
\W_2(\omega) \rightarrow \frac{1-\ee^{-4\omega}}{4\omega}, \qquad
\W_3(\omega) \rightarrow \ee^{-2\omega}\qquad (\rho\downarrow0).
\]
Differentiating the LSTs gives the following results as $\rho\downarrow0$:
\begin{align*}
\E[W_1] &\rightarrow 2,& \E[W_2] &\rightarrow 2, &\E[W_3] &\rightarrow 2, \\
\sd[W_1] &\rightarrow \sqrt{4/3}, &\sd[W_2] &\rightarrow \sqrt{4/3}, &\sd[W_3] &\rightarrow 0.
\end{align*}

\paragraph{Example 2: a two-stage queueing model with customer feedback.}

\begin{figure}[ht]
\begin{center}
\includegraphics[width=0.7\linewidth]{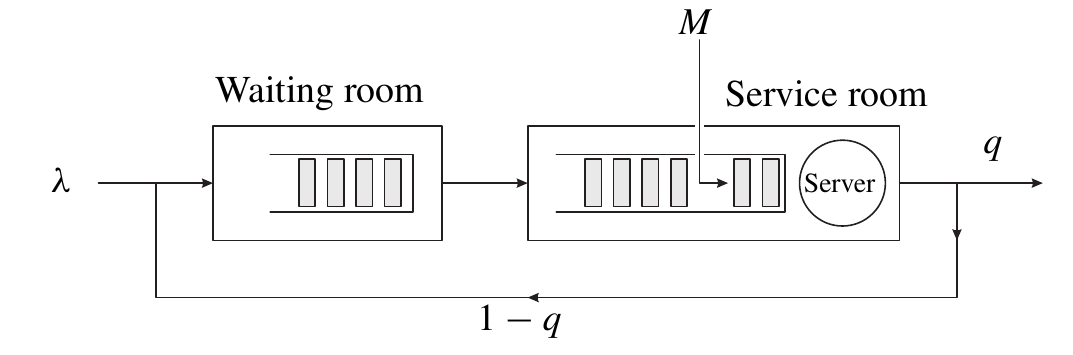}
\end{center}
\caption{The two-stage queueing model with customer feedback, as discussed in Example 2.}
\label{figureExample2}
\end{figure}

This second example is introduced by Tak\'acs \cite{takacsfeedback77}, and extended by Ali and Neuts \cite{alineuts84}. The queueing system under consideration consists of a waiting room, in which customers arrive according to a Poisson process with intensity $\lambda$, and a service room. The customers are all transferred simultaneously to the service room where they receive service in order of arrival. However, at the moment of the transfer to this service room $M$ additional ``overhead customers'' are added to the front of this queue. (In \cite{takacsfeedback77} $M$ is a constant, in \cite{alineuts84} it is a random variable.) Upon service completion, each customer leaves the system with probability $q$, and returns to the waiting room with probability $1-q$. Overhead customers leave the system with probability one after being served. As soon as the last customer in the service room finishes service (and either leaves the system, or returns to the waiting room) all customers present in the waiting room are transferred to the service room, where they will receive service after a new batch of overhead customers has been served, and so on. A schematic representation of this model is depicted in Figure \ref{figureExample2}.

We use the same input parameters as Tak\'acs \cite{takacsfeedback77}: $q=2/3$ and $\lambda/\mu=1/6$, where $1/\mu$ is the mean service time in the service room. This service time is exponentially distributed. The number of overhead customers that are added to the front of the queue is a constant with value $M$. We can model this system in terms of our network with a single, shared server by defining arrival intensities $\lambda_1=\lambda$ and $\lambda_2 = 0$. The service times in stations 1 and 2 are respectively $0$ and exponentially distributed with mean $b_2=1/\mu$. The routing probabilities are $p_{1,2}=1$ and $p_{2,1}=1/3$, the other $p_{i,j}$ are zero. The service times of the overhead customers are also exponentially distributed with parameter $\mu$. Hence, we can model the addition of $M$ overhead customers as a switch-over time which is Erlang-$M$ distributed with parameter $\mu$. The switch-over time between $Q_2$ and $Q_1$ is zero.  Note that, since $b_1=0$, there is no difference between gated and exhaustive service. By differentiation of the waiting time LSTs \eqref{LSTW}, we can obtain explicit expressions for all moments of the waiting-time distributions for this example. The first three moments of the waiting times are given below.
\begin{align*}
&\E[W_1] = \frac{1+M}{2\mu}, &&\E[W_1^2] = \frac{(M+1) (11 M+25)}{27 \mu ^2}, &&\E[W_1^3] =\frac{(M+1) (M (43
   M+223)+310)}{108 \mu ^3},\\
&\E[W_2]=\frac{1+7M}{6\mu}, &&\E[W_2^2]=\frac{(M+1) (37 M+11)}{27 \mu ^2},&&\E[W_2^3]=\frac{(M+1) (M+2) (175
   M+81)}{108 \mu ^3}.
\end{align*}
The results are slightly different from those presented in \cite{takacsfeedback77}, because Tak\'acs also considers the overhead customers in the computations of the waiting times and allows them to return to the waiting room after their service is completed. Modelling this situation would require one minor adaptation in the laws of motion (adding the overhead customers at the beginning of $V_2$) and another adaptation in the waiting time LST (conditioning on the event that a new customer is an overhead customer). These changes are not too difficult but beyond the scope of this paper.

\section{Discussion and further research}\label{discussion}

In this section, we not only elaborate on the developed method and its applicability, but we also discuss possible ways of extending the present study.

\paragraph{Method.}

As mentioned in the introduction, the main complicating factor of the model under consideration is caused by the rerouting of internal customers. This implies that the \emph{total} arrival process at each queue is not Poisson, and not even renewal. Traditional methods to determine waiting-time distributions in each queue are based on the distributional form of Little's Law, which relies on the assumption of Poisson arrivals. Contrary to the distributional form of Little's Law, we explicitly make use of the branching structure to find waiting-time distributions. The main idea is that upon the arrival of a tagged customer $Y$ at time $t$ at $Q_i$ we compute \emph{a priori} the total future service times at each of the queues, for \emph{all} the other customers present in the system at time $t$ that will be served before customer $Y$ enters service at $Q_i$ (see \eqref{bkiomega}). Additionally, we add the total future service requirements of all external arrivals (and their descendants) that will be served before customer $Y$ enters service (see \eqref{rkiomega}). The advantage of this method is that a system no longer needs to satisfy all of the prerequisites required to apply the distributional form of Little's Law (see \cite{keilsonservi90}).

\paragraph{Applicability.}

The novel approach of this paper to find the LST of the waiting-time distribution can also be applied to other types of models with a single server serving multiple queues. Obviously, one can apply it to standard polling models (without customer routing) by simply taking $p_{i,0} = 1$ and $p_{i,j} = 0$ for $j > 0$. However, the developed methodology carries almost directly over to tandem queues \cite{nair,taube}, multi-stage queueing models with parallel queues \cite{katayama}, feedback vacation queues \cite{boxmayechiali97, takine}, symmetric feedback polling systems \cite{takagifeedback,takine}, systems with a waiting room \cite{alineuts84,takacsfeedback77},  closed networks \cite{altman2}, $M/G/1$ queues with permanent and transient customers \cite{boxmacohen91}, networks with permanent and transient customers \cite{armonyyechiali99}, or polling models with arrival rates that depend on the location of the server \cite{boonsmartcustomers2010,smartcustomers}.

\paragraph{Further research.}

Since the model can be described as a multi-type branching process, \emph{explicit closed-form} expressions can be obtained for the waiting-time distributions under heavy-traffic (HT) assumptions. Such expressions are appealing because they give fundamental insight in how the system performance depends on the system parameters, and in particular on the routing probabilities $p_{i,j}$. HT asymptotics can be obtained by combining insights from multi-type branching processes \cite{RvdM_QUESTA}, fluid analyses \cite{olsenvdmei03,olsenvdmei05}, and the heavy-traffic averaging principle by Coffman et al. \cite{coffman95,coffman98}.
The HT analysis is relevant because in practice the proper operation of the system is particularly important when the system is heavily loaded. The HT asymptotics form an excellent basis for the development of approximations for the waiting-time distributions for \emph{arbitrary} loads.
For the \emph{mean} waiting times, preliminary results are obtained in \cite{boonvdmeiwinandsRovingPER2011}.

From a practical perspective, motivated by applications in production systems \cite{boonapplications2011}, an important extension of the model under consideration is a model where customers visit a predetermined, class-specific sequence of queues in a fixed order. In our model one would have to define multiple customer classes, each having their own fixed visit order through the system. The method presented in this paper forms a good basis for this extension.

\section*{Acknowledgements}

The authors are grateful to Ivo Adan and Onno Boxma for providing valuable comments on earlier drafts of the present paper.

\expandafter\ifx\csname urlstyle\endcsname\relax
  \providecommand{\doi}[1]{DOI: #1}\else
  \providecommand{\doi}{DOI: \begingroup \urlstyle{rm}\Url}\fi

\bibliographystyle{abbrvnat}
%\bibliography{roving}

\begin{thebibliography}{35}
\providecommand{\natexlab}[1]{#1}
\providecommand{\url}[1]{\texttt{#1}}
\expandafter\ifx\csname urlstyle\endcsname\relax
  \providecommand{\doi}[1]{doi: #1}\else
  \providecommand{\doi}{doi: \begingroup \urlstyle{rm}\Url}\fi

\bibitem[Ali and Neuts(1984)]{alineuts84}
O.~M.~E. Ali and M.~F. Neuts.
\newblock A service system with two stages of waiting and feedback of
  customers.
\newblock \emph{Journal of Applied Probability}, 21:\penalty0 404--413, 1984.

\bibitem[Altman and Yechiali(1994)]{altman2}
E.~Altman and U.~Yechiali.
\newblock Polling in a closed network.
\newblock \emph{Probability in the Engineering and Informational Sciences},
  8\penalty0 (3):\penalty0 327--343, 1994.

\bibitem[Armony and Yechiali(1999)]{armonyyechiali99}
R.~Armony and U.~Yechiali.
\newblock Polling systems with permanent and transient jobs.
\newblock \emph{Communications in Statistics. Stochastic Models}, 15\penalty0
  (3):\penalty0 395--427, 1999.

\bibitem[Boon et~al.(2010)Boon, van Wijk, Adan, and
  Boxma]{boonsmartcustomers2010}
M.~A.~A. Boon, A.~C.~C. van Wijk, I.~J. B.~F. Adan, and O.~J. Boxma.
\newblock A polling model with smart customers.
\newblock \emph{Queueing Systems}, 66\penalty0 (3):\penalty0 239--274, 2010.

\bibitem[Boon et~al.(2011{\natexlab{a}})Boon, van~der Mei, and
  Winands]{boonapplications2011}
M.~A.~A. Boon, R.~D. van~der Mei, and E.~M.~M. Winands.
\newblock Applications of polling systems.
\newblock \emph{Surveys in Operations Research and Management Science},
  16:\penalty0 67--82, 2011{\natexlab{a}}.

\bibitem[Boon et~al.(2011{\natexlab{b}})Boon, van~der Mei, and
  Winands]{boonvdmeiwinandsRovingPER2011}
M.~A.~A. Boon, R.~D. van~der Mei, and E.~M.~M. Winands.
\newblock Queueing networks with a single shared server: light and heavy
  traffic.
\newblock \emph{SIGMETRICS Performance Evaluation Review}, 39\penalty0
  (2):\penalty0 44--46, 2011{\natexlab{b}}.

\bibitem[Boxma(1994)]{smartcustomers}
O.~J. Boxma.
\newblock Polling systems.
\newblock In K.~Apt, L.~Schrijver, and N.~Temme, editors, \emph{From universal
  morphisms to megabytes: A Baayen space odyssey -- Liber amicorum for P. C.
  Baayen}, pages 215--230. CWI, Amsterdam, 1994.

\bibitem[Boxma and Cohen(1991)]{boxmacohen91}
O.~J. Boxma and J.~W. Cohen.
\newblock The {$M/G/1$} queue with permanent customers.
\newblock \emph{IEEE Journal on Selected Areas in Communications}, 9\penalty0
  (2):\penalty0 179--184, 1991.

\bibitem[Boxma and Yechiali(1997)]{boxmayechiali97}
O.~J. Boxma and U.~Yechiali.
\newblock An {$M/G/1$} queue with multiple types of feedback and gated
  vacations.
\newblock \emph{Journal of Applied Probability}, 34:\penalty0 773--784, 1997.

\bibitem[Boxma et~al.(2009)Boxma, Bruin, and Fralix]{boxmafralixbruin08}
O.~J. Boxma, J.~Bruin, and B.~H. Fralix.
\newblock Waiting times in polling systems with various service disciplines.
\newblock \emph{Performance Evaluation}, 66:\penalty0 621--639, 2009.

\bibitem[Boxma et~al.(2011)Boxma, Kella, and
  Kosi\'nski]{boxmakellakosinski2011}
O.~J. Boxma, O.~Kella, and K.~M. Kosi\'nski.
\newblock Queue lengths and workloads in polling systems.
\newblock \emph{Operations Research Letters}, 39:\penalty0 401--405, 2011.

\bibitem[{Coffman, Jr.} et~al.(1995){Coffman, Jr.}, Puhalskii, and
  Reiman]{coffman95}
E.~G. {Coffman, Jr.}, A.~A. Puhalskii, and M.~I. Reiman.
\newblock Polling systems with zero switchover times: A heavy-traffic averaging
  principle.
\newblock \emph{The Annals of Applied Probability}, 5\penalty0 (3):\penalty0
  681--719, 1995.

\bibitem[{Coffman, Jr.} et~al.(1998){Coffman, Jr.}, Puhalskii, and
  Reiman]{coffman98}
E.~G. {Coffman, Jr.}, A.~A. Puhalskii, and M.~I. Reiman.
\newblock Polling systems in heavy-traffic: A {Bessel} process limit.
\newblock \emph{Mathematics of Operations Research}, 23:\penalty0 257--304,
  1998.

\bibitem[Eisenberg(1972)]{eisenberg72}
M.~Eisenberg.
\newblock Queues with periodic service and changeover time.
\newblock \emph{Operations Research}, 20\penalty0 (2):\penalty0 440--451, 1972.

\bibitem[Foss(1984)]{fossBranching}
S.~Foss.
\newblock Queues with customers of several types.
\newblock In A.~A. Borovkov, editor, \emph{Advances in Probability Theory:
  Limit Theorems and Related Problems}, pages 348--377. Optimization Software,
  1984.

\bibitem[Gong and de~Koster(2008)]{gongdekoster08}
Y.~Gong and R.~de~Koster.
\newblock A polling-based dynamic order picking system for online retailers.
\newblock \emph{IIE Transactions}, 40:\penalty0 1070--1082, 2008.

\bibitem[Grasman et~al.(2008)Grasman, Olsen, and Birge]{Grasman1}
S.~E. Grasman, T.~L. Olsen, and J.~R. Birge.
\newblock Setting basestock levels in multiproduct systems with setups and
  random yield.
\newblock \emph{IIE Transactions}, 40\penalty0 (12):\penalty0 1158--1170, 2008.

\bibitem[Grillo(1990)]{Grillo1}
D.~Grillo.
\newblock Polling mechanism models in communication systems -- some application
  examples.
\newblock In H.~Takagi, editor, \emph{Stochastic Analysis of Computer and
  Communication Systems}, pages 659--699. North-Holland, Amsterdam, 1990.

\bibitem[Katayama(1988)]{katayama}
T.~Katayama.
\newblock A cyclic service tandem queueing model with parallel queues in the
  first stage.
\newblock \emph{Stochastic Models}, 4:\penalty0 421--443, 1988.

\bibitem[Kavitha and Altman(2009)]{kavitha}
V.~Kavitha and E.~Altman.
\newblock Queueing in space: design of message ferry routes in static adhoc
  networks.
\newblock In \emph{Proceedings ITC21}, 2009.

\bibitem[Keilson and Servi(1990)]{keilsonservi90}
J.~Keilson and L.~D. Servi.
\newblock The distributional form of {Little's Law} and the {Fuhrmann-Cooper}
  decomposition.
\newblock \emph{Operations Research Letters}, 9\penalty0 (4):\penalty0
  239--247, 1990.

\bibitem[Levy and Sidi(1990)]{levy1}
H.~Levy and M.~Sidi.
\newblock Polling systems: applications, modeling, and optimization.
\newblock \emph{IEEE Transactions on Communications}, 38:\penalty0 1750--1760,
  1990.

\bibitem[Nair(1971)]{nair}
S.~S. Nair.
\newblock A single server tandem queue.
\newblock \emph{Journal of Applied Probability}, 8\penalty0 (1):\penalty0
  95--109, 1971.

\bibitem[Olsen and van~der Mei(2003)]{olsenvdmei03}
T.~L. Olsen and R.~D. van~der Mei.
\newblock Polling systems with periodic server routeing in heavy traffic:
  distribution of the delay.
\newblock \emph{Journal of Applied Probability}, 40:\penalty0 305--326, 2003.

\bibitem[Olsen and van~der Mei(2005)]{olsenvdmei05}
T.~L. Olsen and R.~D. van~der Mei.
\newblock Periodic polling systems in heavy-traffic: renewal arrivals.
\newblock \emph{Operations Research Letters}, 33:\penalty0 17--25, 2005.

\bibitem[Resing(1993)]{resing93}
J.~A.~C. Resing.
\newblock Polling systems and multitype branching processes.
\newblock \emph{Queueing Systems}, 13:\penalty0 409--426, 1993.

\bibitem[Sarkar and Zangwill(1992)]{sarkar}
D.~Sarkar and W.~I. Zangwill.
\newblock File and work transfers in cyclic queue systems.
\newblock \emph{Management Science}, 38\penalty0 (10):\penalty0 1510--1523,
  1992.

\bibitem[Sidi and Levy(1990)]{sidi1}
M.~Sidi and H.~Levy.
\newblock Customer routing in polling systems.
\newblock In P.~King, I.~Mitrani, and R.~Pooley, editors, \emph{Proceedings
  Performance '90}, pages 319--331. North-Holland, Amsterdam, 1990.

\bibitem[Sidi et~al.(1992)Sidi, Levy, and Fuhrmann]{sidi2}
M.~Sidi, H.~Levy, and S.~W. Fuhrmann.
\newblock A queueing network with a single cyclically roving server.
\newblock \emph{Queueing Systems}, 11:\penalty0 121--144, 1992.

\bibitem[Tak\'acs(1977)]{takacsfeedback77}
L.~Tak\'acs.
\newblock A queuing model with feedback.
\newblock \emph{Revue fran\c{c}aise d'automatique, d'informatique et de
  recherche op\'erationnelle. Recherche op\'erationnelle}, 11\penalty0
  (4):\penalty0 345--354, 1977.

\bibitem[Takagi(1987)]{takagifeedback}
H.~Takagi.
\newblock Analysis and applications of a multiqueue cyclic service system with
  feedback.
\newblock \emph{IEEE Transactions on Communications - TCOM}, 35\penalty0
  (2):\penalty0 248--250, 1987.

\bibitem[Takagi(2000)]{takagi3}
H.~Takagi.
\newblock Analysis and application of polling models.
\newblock In G.~Haring, C.~Lindemann, and M.~Reiser, editors, \emph{Performance
  Evaluation: Origins and Directions}, volume 1769 of \emph{Lecture Notes in
  Computer Science}, pages 424--442. Springer Verlag, Berlin, 2000.

\bibitem[Takine et~al.(1991)Takine, Takagi, and Hasegawa]{takine}
T.~Takine, H.~Takagi, and T.~Hasegawa.
\newblock Sojourn times in vacation and polling systems with {Bernoulli}
  feedback.
\newblock \emph{Journal of Applied Probability}, 28\penalty0 (2):\penalty0
  422--432, 1991.

\bibitem[Taube-Netto(1977)]{taube}
M.~Taube-Netto.
\newblock Two queues in tandem attended by a single server.
\newblock \emph{Operations Research}, 25\penalty0 (1):\penalty0 140--147, 1977.

\bibitem[Van~der Mei(2007)]{RvdM_QUESTA}
R.~D. Van~der Mei.
\newblock Towards a unifying theory on branching-type polling models in heavy
  traffic.
\newblock \emph{Queueing Systems}, 57:\penalty0 29--46, 2007.

\end{thebibliography}

\end{document}